\documentclass[12pt]{article}

\usepackage{a4,amsmath,amsfonts, amssymb}

\begin{document}
\bibliographystyle{plain}

\def \outlineby #1#2#3{\vbox{\hrule\hbox{\vrule\kern #1%
\vbox{\kern #2 #3\kern #2}\kern #1\vrule}\hrule}}%
\def \endbox {\outlineby{4pt}{4pt}{}}%

\newenvironment{proof}
{\noindent{\bf Proof\ }}{{\hfill \endbox
}\par\vskip2\parsep}

\hfuzz5pt
\setlength{\parindent}{0pt}

\newcommand{\vt}{\widetilde{V}(m)}
\newcommand{\wt}{\widetilde{W}(m)}
\newcommand{\ep}{\mathbb E}
\newcommand{\I}{\mathbb I}
\newcommand{\pr}{\mathbb P}
\newcommand{\re}{\mathbb R}
\newcommand{\tends}{\rightarrow \infty}
\newcommand{\var}{\mbox{Var}}
\newcommand{\cov}{\mbox{Cov}}
\newcommand{\nm}[1]{\langle #1 \rangle }
\newcommand{\rex}{\re \cup \{ \infty \} }
\newcommand{\ene}{\varepsilon}
\newcommand{\tr}{\mbox{tr}}
\newcommand{\rat}[2]{\log \left( \frac{ #1}{#2} \right)}
\newcommand{\supp}{\mbox{{\em{supp}}}}
\newcommand{\ptnfn}[1]{Z_{#1}}
\newcommand{\ptn}{\ptnfn{c}}
\newcommand{\gibb}[1]{g_{#1,c}}
\newcommand{\gib}{\gibb{n}}
\newcommand{\vc}[1]{{\mathbf{#1}}}

\newtheorem{theorem}{Theorem}[section]
\newtheorem{lemma}[theorem]{Lemma}
\newtheorem{proposition}[theorem]{Proposition}
\newtheorem{corollary}[theorem]{Corollary}
\newtheorem{conjecture}[theorem]{Conjecture}
\newtheorem{definition}[theorem]{Definition}
\newtheorem{example}[theorem]{Example}
\newtheorem{condition}{Condition}
\title{Entropy and a generalisation\\ of ``Poincare's Observation''}
\author{Oliver Johnson\thanks{Statistical Laboratory, 
CMS, Wilberforce Road, Cambridge, CB3 0WB, UK. Email: {\tt otj1000@cam.ac.uk}. 
Fax: +44 1223 337956}}
\date{\today}
\maketitle
\begin{abstract} \noindent
Consider a sphere of radius $\sqrt{n}$ in $n$ dimensions, and
consider $\vc{X}$, a random variable uniformly distributed on its surface.
Poincar\'{e}'s Observation states that for large $n$, 
the distribution of the first $k$ coordinates of $\vc{X}$
is close in total variation distance to the standard normal $N(\vc{0},I_k)$. 
In this paper, we consider a larger family of manifolds, and $\vc{X}$
taking a more general distribution on the surfaces. We establish a bound
in the stronger Kullback--Leibler sense of relative entropy, and discuss
its sharpness, providing a necessary condition for convergence in this sense.
We show how our results imply the equivalence of ensembles for a wider 
class of test functions than is standard. We also 
deduce results of de Finetti type, concerning a generalisation of the idea of
orthogonal invariance.
\renewcommand{\thefootnote}{}
\renewcommand{\thefootnote}{\arabic{footnote}}
\setcounter{footnote}{0}
\end{abstract}
\section{Notation and Definitions}
Diaconis  and  Freedman \cite{diaconis}  consider  $\vc{X}$, a  random
variable uniformly  distributed on the  surface of a sphere  of radius
$\sqrt{n}$  in $n$  dimensions.  They show  that  the distribution  of
$\pi_{n,k} \vc{X}$, the first $k$  coordinates of $\vc{X}$, is close in
total variation distance to a  Gaussian for large $n$.  They indicate
a  natural connection  between  this problem  and  the equivalence  of
ensembles for  the free Hamiltonian,  where a sphere corresponds  to a
surface of  constant kinetic  energy. In this  paper, we  consider a
more  general  family  of  manifolds defined  by  symmetric,  additive
Hamiltonians.  Unlike Diaconis  and Freedman we shall not assume
a uniform distribution on  these surfaces.  We  prove convergence in
the stronger Kullback--Leibler sense of relative entropy distance.

First, we introduce notation. Given measurable $f:\re \rightarrow \rex$ 
we define $F$ to be the set on which $f$ is finite,  and let
$R_n = R_n(f)$ be given by $R_n(\vc{x}) = \sum_{i=1}^{n} f(x_i),$
We define the surface ${\mathcal{S}}_n(t) = \{ \vc{x}: R_n(\vc{x}) = nt \}
\subseteq F^n.$
Given $f$ and $c>0$, define $\gib$ for the density of a Gibbs distribution:
$ \gib(\vc{x}) =       \exp(-cR_n(\vc{x}))/ \ptn^n. $
We assume that $f$ has the property that $\ptn = \int \exp(-cf(x)) 
dx$ is finite for all $c>0$.
Define the projection $\pi_{n,k}: \re^n \rightarrow \re^k$ restricting to 
the first $k$ coordinates: 
$ \pi_{n,k} (x_1, x_2, \ldots, x_n) = (x_1, x_2, \ldots, x_k)$.
For any probability density $q$ on $\re^n$
write the energy $\ene(q) = \int q(\vc{x}) R_n(\vc{x}) d\vc{x}$ and 
the entropy $h(q) = \int - q(\vc{x}) \log q(\vc{x}) d\vc{x}$.

For a given $f$, we know that Gibbs densities
maximise the entropy for a given energy. (This characterisation of Gibbs 
measures via a variational principle is discussed
in Chapter 15 of Georgii \cite{georgii2}). This is done by considering the
Kullback--Leibler distance $D( p\| \gibb{n} )$, 
where $D(f\| g) = \int f(\vc{x}) \log 
(f(\vc{x})/g(\vc{x})) d\vc{x}$, for $f$ and $g$ probability 
densities on $\re^k$.
\begin{example} \label{eg:disns} Two cases in particular are significant here:
\begin{enumerate}
\item{If $f(x) = x^2$, then $\ptn = \sqrt{\pi/c}$, and 
$\gib(\vc{x}) = \prod_{i=1}^n \gibb{1}(x_i)$, where $\gibb{1}$ is a 
$N(0, (2c)^{-1})$ density, 
the density which maximises entropy subject to a variance 
constraint. \label{eg:norm} }
\item{If $f(x) = x$ for $x \geq 0$ and $f(x) = \infty$ for $x < 0$, then
$\ptn = 1/c$, and $\gib(\vc{x}) = \prod_{i=1}^n \gibb{1}(x_i)$,
where $\gibb{1}$ is an Exp($c$) density, the density which maximises 
entropy on the positive half-line subject to a mean constraint.
\label{eg:exp} }
\end{enumerate} \end{example}
Now we can state Diaconis and Freedman's results on the normal and 
exponential case in the form:
\begin{theorem} \label{thm:diaconis}
For $f$ as in Example \ref{eg:disns}.\ref{eg:norm} or 
\ref{eg:disns}.\ref{eg:exp} and for a given $t \in (0, \infty)$, 
take $c$ such that $\ene(\gibb{1})= t$ and let 
$\vc{X}$ be distributed uniformly on the 
surface ${\mathcal{S}}_n (t) \subseteq F^n$.
Writing $p_{n,k,t}(\vc{y})$ for  
the density of $\pi_{n,k} \vc{X}$:

For $f(x) = x^2$, as in Example \ref{eg:disns}.\ref{eg:norm},
$d_{TV}(p_{n,k,t}, \gibb{k}) \leq 2(k+3)/(n-k-3)$.

For $f(x) = x$, as in Example \ref{eg:disns}.\ref{eg:exp},
$d_{TV}(p_{n,k,t} \vc{X}, \gibb{k}) \leq 2(k+1)/(n-k-1)$.
\end{theorem} 
The following theorem is the main result of this paper. It involves
a class $\mathcal{F}$ of functions (we postpone the precise definition to 
Definition \ref{def:f}, but roughly speaking we want $f$ to be strictly
increasing and well-behaved at zero). Let $h_{n,t}$ be the density 
corresponding to uniform distribution on ${\mathcal{S}}_n(t)$.
\begin{theorem} \label{thm:convind}
Assuming $f \in \mathcal{F}$, for any $t \in (0, \infty)$ 
take $c$ such that $\ene(\gibb{1})= t$ and consider $\vc{X}$, a random variable
with  density $p$ on 
${\mathcal{S}}_n(t)$. Writing $p_{n,k,t}(\vc{y})$ for  
the density of $\pi_{n,k} \vc{X}$, for some
constant $C = C(f)$:
$$ D(p_{n,k,t} \| \gibb{k} ) \leq D(p \| h_{n,t}) + \rat{n}{n-k}
+ \frac{2}{\sqrt{n}/C - 1}.$$\end{theorem}
A major motivation for this work, beyond the intrinsic interest of 
generalising Diaconis and Freedman's work \cite{diaconis}, 
comes from the question of so-called `equivalence of ensembles'.
Given a Hamiltonian $H(\vc{x})$, for $\vc{x} \in \re^n$  we  can consider  the
microcanonical  ensemble  (uniform measure with density $h_{n,t}$ on
the surface $\{\vc{x}: H(\vc{x})=nt \}$) and grand canonical ensemble
(Gibbs measure with density $g_{n,c}$   
proportional to $\exp(-cH(\vc{x}))$ for $\vc{x} \in \re^n$). The principle of
equivalence  of  ensembles   suggests that in some sense $h_{n,t}$ and 
$g_{n,c}$ are close together, that is 
\begin{equation} \label{eq:eqens}
\int p(\vc{x}) h_{n,t}(d\vc{x}) \simeq \int p(\vc{x}) g_{n,c}(\vc{x})
\end{equation}
for some class of test functions $p$.
Convergence in  total variation  distance (as established  by Diaconis
and  Freedman) implies  that Equation  (\ref{eq:eqens}) holds  for $p$
bounded  and continuous and depending only on $k = o(n)$ coordinates. 
Theorem \ref{thm:convind} implies this for a wider  class of test functions,
including the Hamiltonian itself.
\begin{corollary}
Assuming $f \in \mathcal{F}$, for any $t \in (0, \infty)$ 
$$ \lim_{n \tends} \left( \int p(\vc{x}) h_{n,t}(d\vc{x}) -
\int p(\vc{x}) g_{n,c}(\vc{x}) \right) = 0,$$
for $p$ depending only on $k$ coordinates, if
$p(x_1,\ldots, x_k)$ is bounded above by a multiple of $1+ \sum_{i=1}^k 
f(x_i)$.
Here, $k$ need not be fixed, and can grow as $o(n)$.
\end{corollary}
\begin{proof} This is a consequence of Theorem \ref{thm:convind} in
conjunction with Lemma 3.1 of Csisz\'{a}r \cite{csiszar5}. The latter states
that $D(u_n \| v) \rightarrow 0$ implies $\int p(x) u_n(x) dx \rightarrow
\int p(x) v(x) dx$, for any $p$ such that $\int \exp(tp(x)) v(x) dx < \infty$
for small $|t|$. This integral is seen to be finite by definition of
the partition function. \end{proof}

A second application is to results of de Finetti type, as described in
\cite{diaconis}. For example,
define an infinite sequence $X_1, X_2, \ldots   $  of random  variables 
to be orthogonally  invariant if for any $n$, 
the law of $(X_1, \ldots, X_n)$ is invariant 
under orthogonal transformations of $\re^n$. Schoenberg \cite{schoenberg} 
states that all orthogonally invariant distributions are mixtures of normals. 
Diaconis and Freedman \cite{diaconis} prove this by showing that the 
first $k$ of 
$n$ orthogonally invariant variables are within $2k/n$ of a mixture of normals,
so a passage to the limit provides the infinite result.

In a similar way we can consider $f$-invariant measures; that is, sequences
of random variables such that for any $n$, the law of $(X_1, \ldots, X_n)$ is 
invariant under continuous transformations of $\re^n$ which preserve 
$R_n(x_1, \ldots, x_n)$. We show:
\begin{corollary}
For $f \in {\mathcal F}$, the only $f$-invariant measures $P$ are
mixtures of Gibbs measures. That is, if we write $G_{\infty,c}$ for
the distribution of $Y_1, Y_2, \ldots$, where $Y_i$ are independent 
with density $\gibb{1}$, there exists a measure $\lambda$, valued
on $c >0$, such that
$$ P= \int G_{\infty,c} \lambda(c) dc.$$
\end{corollary}
\begin{proof} First we consider finite subsequences of $X_1, \ldots$.
If measure $P_n$ is invariant under $R_n$-preserving transformations then
(as in \cite{diaconis}), $P_n$ is constant on each manifold 
${\cal S}_n(t)$, so:
$$ P_n = \int h_{n,t} \mu_n(t) dt,$$
for some measure probability measure
$\mu_n$. Now projecting down to the first $k$ coordinates:
$$ Q_k = \pi_{n,k} P_n = \int (\pi_{n,k} h_{n,t}) \mu_n(t) dt.$$
Writing $c(t)$ for the unique $c$ such that $\ene(g_{1,c}) =t$, and
defining $R_k = \int g_{k,c(t)} \mu_n(t) dt$, then:
$$ \| R_k - Q_k \|_{TV} \leq \int \mu_n(t) \| (\pi_{n,k} h_{n,t}) - g_{k,c(t)}
\| dt \leq \sqrt{\frac{2k}{n-k}},$$
by Theorem \ref{thm:convind}.

Duplicating Diaconis and Freedman's tightness argument we can show that 
as $n \rightarrow \infty$, $\mu_n$ must have a convergent subsequence,
with limit $\mu$. Mapping from $t$ to $c(t)$, we find a $\lambda$ with
the required properties.
\end{proof}

The fact that projecting a uniform distribution on a
 sphere approximately gives a Gaussian
is often referred to as Poincar\'{e}'s Observation (see for
example \cite{ledoux}). However, in Section 6 of their paper, Diaconis and
Freedman suggest that this attribution is wrong, and  
that the earliest reference to it in the probability 
literature comes in the work of Borel. Nonetheless, it appears that the 
observation is even older than this, and can be traced back to Mehler 
\cite{mehler}.

Csisz\'{a}r \cite{csiszar4} uses entropy-theoretic methods to consider 
the distribution of a random variable $X_1$, conditional on the vector
$\vc{X} = (X_1, \ldots, X_n) \in A$ for some set $A$. However, these results
rely on $A$ being a `thick set' of positive measure, whereas in the present
paper we consider the so-called `thin shell' case. Dembo and Zeitouni 
\cite{dembo3} extended Csisz\'{a}r's results to the distribution
of the first $k$ coordinates,
where $k = k(n)$ can vary with $n$, and (as in this paper) discovered that a 
sufficient condition for convergence is that $k(n) = o(n)$. Schroeder 
\cite{schroeder} generalised this result to Markov processes, and Comets
and Zeitouni \cite{comets} even to mean--field perturbations of Markov 
processes. However, all these papers use the assumption that $A$ has non-empty
interior. This corresponds to a weak limit theorem, bounding the probability
of a particular set, whereas our methods require a local limit theorem, 
bounding densities.
\section{Class of Functions Considered}
It is natural to ask for the widest possible class of functions
$f$ such that results such as Theorem \ref{thm:diaconis} hold.
For technical reasons, we will need to control  $w_n$, the density of
$R_n(\vc{X})$,  when $\vc{X}$ has Gibbs
density $\gib$.   Specifically we need an
upper bound on  $\log(w_{n-k}(s)/w_n(t))$ for certain $s,t$.  Theorem
\ref{thm:diaconis} holds because  $w_n$  is known exactly (in
Example   \ref{eg:disns}.\ref{eg:norm} $R_n(\vc{X})$  is  the sum   of
squares of normals, with the $\chi^2_n$ distribution, and in Example
\ref{eg:disns}.\ref{eg:exp} $R_n(\vc{X})$  is   the sum of exponentials,   
with the $\Gamma(n,c)$ distribution).  A later paper by Borovkov 
\cite{borovkov} extends Diaconis and Freedman's work to the case 
$f(x) = |x|^p$, again using exact calculations of the density. 
Whilst his bounds are tighter, the method will not extend to the general
case -- it gives no information about Hamiltonians of 
the form $f(x) =x^2+\epsilon x^4$.

However, we  can make progress in  other cases too.
The key  observation  is that $X_i$  are  IID (with   marginal density
$\gibb{1}$), so $R_n(\vc{X})  = \sum_{i=1}^{n} f(X_i)$ is  a sum of $n$ IID
random variables. Using a local version of the Central Limit Theorem, 
we will  be able to show that
the  densities are sufficiently close  in supremum norm to a Gaussian
density for our result to go through.
We obtain a proposition reminiscent of 
Equation (2.9) of \cite{diaconis}.
\begin{proposition} \label{prop:df}
Consider $\vc{X}$ with density $\gib$ and let $Y_j =f(X_j)$, with 
$\ep Y_j = \mu$, $\var(Y_j) = \sigma^2$, and  
characteristic function $\phi(u) = \ep \exp(iuY_j)$.
Assume that $I = \int |\phi(u)|^r du$ is finite for some $r \geq 1$ and  
$m = \ep |Y - \mu|^3$ is finite. 
If $w_k$ is the density of $R_k(\vc{X}) = \sum_{i=1}^k Y_i$ then there
exists a constant $C = C(Y)$ such that for $n-k \geq r$:
$$ \rat{w_{n-k}(z)}{w_n(n \mu)} \leq \rat{n}{n-k} + 
\frac{2 }{\sqrt{n}/C - 1} \mbox{ for all $z \in \re$.}$$
\end{proposition}
\begin{proof}
The local limit theorem tells us that there exists a constant $C(Y)$ such
that if $f_n$ is the density of $\sum (Y_i - \mu)/\sqrt{n\sigma^2}$ then
for $n \geq r$:
$$ \sup_x \left| f_n(x) - \frac{1}{\sqrt{2\pi}} \exp(-x^2/2) \right| 
\leq \frac{C(Y)}{\sqrt{2 \pi n}}.$$
A careful reading of (for example) Section 46 of Gnedenko and Kolmogorov
\cite{gnedenko}
shows that the dependence of $C(Y)$ on $Y$ comes through $m, \sigma^2, I$
and $\nu = \sup_{t > \sigma^2/m} |\phi(t)| < 1$.

Rescaling this, we know that for any $z$:
$$ w_{n-k}(z) \leq \frac{\sqrt{n-k} + C}
{(n-k)\sqrt{2 \pi \sigma^2}} \mbox{ and }
 w_{n}(n \mu) \geq \frac{\sqrt{n} - C}{n \sqrt{2 \pi \sigma^2 }} .$$
Taking the ratio of these terms we deduce the result, since:
$$ \rat{\sqrt{n-k}+C}{\sqrt{n}-C} \leq 
\log \left( 1+ \frac{2C}{\sqrt{n}-C} \right).$$
\end{proof}
We can now describe the class of surfaces that our Theorem will cover. 
These conditions on function $f$ are chosen so that they imply that the 
local limit theorem will hold in Proposition \ref{prop:df}. We confirm this
in Lemmas \ref{lem:meanbd} and \ref{lem:intbd} below. Roughly
speaking, we require the function $f$ to grow faster than a linear function
on most of the domain, and to be well-behaved at $0$.
\begin{definition} \label{def:f}
Define $\mathcal{F}$ to be the class of functions $f:\re \rightarrow 
\rex$ with: \begin{enumerate}
\item{$f(0) = 0$, $f(x)$ right continuous at $0$.}
\item{$f(x)$ is differentiable for $x>0$ with $f'(x) > 0$ and:
\begin{enumerate}
\item{There exist $a_1$,$a_2$ such that $f'(x) \geq a_1 > 0$ for 
$x \in (a_2,\infty)$. \label{item:taildiff} }
\item{\label{item:enddiff} There exist $1 < q < 2$ and $a_3 > 0$ such that:
 $$ \lim \inf_{x \searrow 0} \frac{f'(x)^q}{f(x)} \geq a_3.$$}
\end{enumerate}}
\item{Either (a) $f(x) \equiv \infty$ for $x<0$  or (b) $f(x) = f(-x)$ for all 
$x$. As before, we write $F$ for the interval on which $f$ is finite.
\label{item:cont} } 
\end{enumerate}
\end{definition}
Notice that for any $p \geq 1$, picking $f(x) = x^p$ on $x \geq 0$ and
infinity elsewhere, or $f(x) = |x|^p$  everywhere ensure that $f \in
\mathcal{F}$. Thus the cases considered by Diaconis and Freedman 
\cite{diaconis} and Borovkov \cite{borovkov} are included in our theorems.

Conditions 1 and 2 ensure that $\ptn$ is finite and non-zero for $c \in (0,
\infty)$, since then
$\int \exp(-cf(x)) dx \leq a_2 + \int_{a_2}^{\infty} f'(x) \exp(-cf(x))/a_1
\leq a_2 + \exp(-cf(a_2))/ca_1$, and continuity provides boundedness
away from zero.
\begin{lemma} Assume $f \in {\mathcal{F}}$. If $X$ has density $\gibb{1}(x)
= \exp(-cf(x))/\ptn$ and $Y = f(X)$, then $m = \ep |Y- \ep Y|^3$ is finite.
\label{lem:meanbd} \end{lemma}
\begin{proof}
In case \ref{item:cont}(a) and \ref{item:cont}(b) of Definition \ref{def:f}
$f(X)$ will have the same density, hence we need only
consider case \ref{item:cont}(a). Further, note that
in this case, since $Y \geq 0$, $m = \ep (Y - \mu)^3$, the 3rd cumulant.
Since the moment generating function of $Y$ is $M(t) = Z_{c-t}/Z_c$, the
cumulant generating function $\log M(t) = \log Z _{c-t}- \log Z_c$, and
the result follows.
\end{proof}
\begin{lemma} If $f \in {\mathcal{F}}$, then if $X$ has density 
$\gibb{1}(x)=\exp(-cf(x))/\ptn$ and $Y = f(X)$, with characteristic 
function $\phi$, then there exists $r$ such that $I = \int |\phi(u)|^r du $ 
is finite. \label{lem:intbd} \end{lemma}
\begin{proof} The density $g(y)$ of $Y$ is $Z_c^{-1} \exp(-cy) /f'(f^{-1}(y))$.
As described in for example Theorem 74 of \cite{titchmarsh}, if 
$g \in L^p$, for some $1 < p \leq 2$ then $I$ is finite for $r \geq p/(p-1)$. 
By Condition \ref{item:enddiff} of Definition \ref{def:f}, for $y$ close to 
zero $g(y) \leq \mbox{const.} \exp(-cy) y^{1/q}$,
so picking $p > q$ ensures the finiteness of the integral at 0. Condition
\ref{item:taildiff} gives us control as $y \tends$.
 \end{proof}
\section{Nesting and Projection of Surfaces} \label{sec:nesting}
Recall that we want to project a density $p$ from a manifold 
${\mathcal{S}}_n(t)$ onto its first $k$ coordinates. 
Given $p$, we will first create $(\theta p)$, a density on $F^n$, by 
placing weighted 
copies of $p$ on each manifold ${\mathcal{S}}_n(u)$. One has a 
free choice of how to weight the individual manifolds. The picture is 
that of fitting together an infinite set of `Russian 
dolls', each of different sizes, but each with the same pattern of densities 
on their surface. We then consider the projection of $p$ by considering the 
projection of $(\theta p)$ conditional on being on ${\mathcal{S}}_n(t)$.

By a suitable choice of weighting, we can make the projection well-behaved.
For example, in the case of the uniform density on the sphere, we weight
concentric spheres by the $\chi^2$ distribution to produce the normal density
on $\re^n$. 

We need to develop a coordinate system that allows us to describe a point
$\vc{x}$ in space by giving its distance from the origin $\vc{0}$ and 
the point where the line from $\vc{0}$ to $\vc{x}$ crosses 
${\mathcal{S}}_n(t)$. In   the case of the sphere $f(x)=x^2$, this 
corresponds to transforming between rectangular and polar coordinates.
We shall require one further technical lemma, not proved here:
\begin{lemma} For any $t \in (0,\infty)$, the equation $ 
\ene(\gibb{1}) =t $ has a unique solution as an equation in $c$, if
$f(0) =0$, $f$ is right continuous at zero, increasing on $x > 0$, 
$f(x) \tends$ as $x \tends$, and $f(x) = \infty$ for $x <0$. 
\end{lemma}
More formally, we give a $C^1$-foliation of 
$F^n \backslash \{ \vc{0} \}$ by compact $(n-1)$-dimensional
manifolds ${\mathcal{S}}_n(u)$, $u >0 $. Each manifold is 
endowed with a standard Riemannian metric and
diffeomorphic, by a central projection, 
to either a unit sphere centered
at the origin or to its intersection with a non-negative
orthant. 
\begin{proposition}
For given $t$, under Conditions  1 and 2 of Definition \ref{def:f}
there exists a bijection $\Phi_{t}$ between 
$\{ x: x>0 \} \times {\mathcal{S}}_n(t)$ and
$F^n \setminus \{ \vc{0} \}$. The bijection $\Phi_{t}$ has a Jacobian 
$A_{n,t}$ which is positive everywhere.
\end{proposition}
\begin{proof} For given $u$ and for any $\vc{x} \in F^n \setminus 
\{ \vc{0}\}$, under Conditions  1 and 2 of Definition \ref{def:f}, 
the equation $R_n(k\vc{x}) =u$ has a unique solution in $k >0$, by the 
Intermediate Value Theorem.
Hence given $(r,\vc{s})$, where $r > 0$ and $\vc{s} \in 
{\mathcal{S}}_n(t)$, we can find a unique $\vc{x} = k \vc{s}$, where
$k$ is chosen such that $R_n(k\vc{s})= nr$.
Conversely, given $\vc{x}$, we can define the central projection
$S_{n,t}(\vc{x}) = k\vc{x}$, where $k$
is chosen such that $R_n(k\vc{x}) = nt$, which is equivalent to saying that
$k \vc{x} \in {\mathcal{S}}_n(t)$. We take the pair $(r,\vc{s})
= (R_n(\vc{x}), S_{n,t}(\vc{x}))$.

This foliation induces a coordinate system on $F^n \backslash \{ 
\vc{0} \}$ via the bijection $\Phi_{t}$. 
Let $\sigma_{n,t}$ denote the induced Riemannian volume on 
${\mathcal{S}}_n(t)$. 
The Jacobian
$A_{n,t}(\vc{x}) = A_{n,t}(t,\vc{s})$ is determined since
for any measurable function $f$:
\begin{eqnarray*} \int_{\re_+ \times {\mathcal{S}}_n(t)} f(t,\vc{s}) 
\Phi_{t}(dt,d\vc{s}) & = & 
\int_{\re^n} f(R_n(\vc{x}),S_{n,t}(\vc{x})) d\vc{x} \\
& = & \int_{\re_+ \times {\mathcal{S}}_n(t)} f(t,\vc{s}) 
A_{n,t}(t,\vc{s}) dt \sigma_{n,t}(d\vc{s}).\end{eqnarray*}

The positivity of $f'$ ensures the positivity of $A_{n,t}$, since 
the local structure of the above foliation (and the induced map 
$\Phi_{t}$) may be also described as follows. Given a point 
$\vc{x}$, such that $\Phi_{t}(\vc{x})=  
(R_n(\vc{x}), S_n(\vc{x}))$, for small $\delta \vc{x}$, we calculate
$\Phi_{t}(\vc{x} +\vc{\delta x})$. Now
$R_n(\vc{x} +\vc{\delta x}) = R_n(\vc{x}) + \sum_i \delta x_i f'(x_i)
+ O(\delta x^2)$, which is not identically equal to $R_n(\vc{x})$.
Similarly, if $S_n(\vc{x}) = k\vc{x} \in {\mathcal{S}}_n(t)$, then
$S_n(\vc{x} + \vc{\delta x})$ is $k(1-\epsilon)(\vc{x} +\vc{\delta x})$, 
where $\epsilon$ is chosen such that this lies on ${\mathcal{S}}_n(t)$. 
The choice of $\epsilon$ that achieves this is $\epsilon = 
(\sum f'(x_i) \delta x_i) / (\sum f'(x_i) x_i)$, which again ensures that
$ S_n(\vc{x} + \vc{\delta x})$ is not identically equal to $S_n(\vc{x})$.
 \end{proof}
Now having developed our coordinate system, we can describe the map
$\theta$ which takes a density on ${\mathcal{S}}_n(t)$ and gives a density
on $F^n$. The motivation for this definition is that it gives an 
isometry between densities (see Lemma \ref{lem:isom}).
\begin{definition} 
Given a  probability  density  $p$  on  ${\mathcal{S}}_n(t)$, we   can
define the product density $(w_n \times p)$ on $\re_+ \times 
{\mathcal{S}}_n(t)$, where $w_n$ is the density of $R_n(\vc{X})$ 
when $\vc{X}$ has density $g_{n,c}$. We can thus define the density 
$(\theta p)$ induced by $\Phi_{t}$ on $F^n \setminus \{\vc{0}\}$, since
transforming to Cartesian coordinates, we know that:
$$ (\theta p)(\vc{x})  =  \frac{w_n(R_n(\vc{x}))  p(S_{n,t}(\vc{x}))}{   
A_{n,t}(R_n(\vc{x}),S_{n,t}(\vc{x}))}, \mbox{ for $\vc{x} \neq 
\vc{0}$,}$$
which is by construction a probability density on $F^n \setminus \{\vc{0}\}$.
\end{definition}
The one significant difference, as mentioned by Borovkov, is
that we will no longer consider the uniform distribution on the surface
but rather consider $h_{n,t}(\vc{v}) = A_{n,t}(nt,\vc{v})/ 
\left( \int A_{n,t}(nt,\vc{s}) 
\sigma_{n,t}(d\vc{s}) \right)$ for $\vc{v} \in {\mathcal{S}}_n(t)$.
\begin{lemma} \label{lem:isom} 
With the definitions above:
$D(\theta p \| \gib) = D(p \| h_{n,t})$.
\end{lemma}
\begin{proof} 
Note that $w_n$ is characterized by:
\begin{eqnarray*} \int_{0}^r w_n(t) dt 
& = & \int g_{n,c}(\vc{x}) \I(R_n(\vc{x}) \leq r) d\vc{x}\\
& = & \int \exp(-ct) Z_c^{-n} \I(t \leq r) \left( \int A_{n,t}(nt,\vc{s}) 
\sigma_{n,t}(d\vc{s}) \right) dt. 
\end{eqnarray*}
We deduce that for any $\vc{u} \in F^n$, 
$g_{n,c}(\vc{u}) = w_n(R_n(\vc{u})) / \left( \int
A_{n,t}(nt,\vc{s}) \sigma_{n,t}(d\vc{s}) \right)$.

Hence by definition,
$(\theta h_{n,t})(\vc{u}) = w_n(R_n(\vc{u})) / \left( \int
A_{n,t}(nt,\vc{s}) \sigma_{n,t}(d\vc{s}) \right) = g_{n,c}(\vc{u})$.
This means that 
\begin{eqnarray*}
D(\theta p \| \gib)  & = & D(\theta p  \| \theta h_{n,t}) 
=  \int \theta p(\vc{x}) 
\log \left( \frac{ p(S_{n,t}(\vc{x}))}{h_{n,t}(S_{n,t}(\vc{x}))} 
\right) d\vc{x} \\
& = & \int p(\vc{s}) 
\log \left( \frac{ p(\vc{s})}{h_{n,t}(\vc{s})} \right) d\vc{s} 
\end{eqnarray*}
as required.
\end{proof}
\begin{definition} \label{def:density}
Given $t > 0$ and a density $p$ on ${\mathcal{S}}_n(t)$, 
write $p_{n,k,t}(\vc{y})$ for the density of $p$ projected by 
$\pi_{n,k}$ onto $\re^k$.
\end{definition}
The key observation is that
$   p_{n,k,t}(\vc{y})   = q_{n,k}(\vc{y},nt)/w_n(nt),$  
which is the projection of $(\theta p)$ on $\re^k$ conditioned on $R_n =nt$.
Here  $q_{n,k}(\vc{y},nt)$  may be treated as
the joint density of  $\pi_{n,k}\vc{X}=\vc{y}$ and $R_n(\vc{X})  =nt$,
where $\vc{X}$ is a random point of $F^n$ with density $\theta p$.
We thus establish the principal result of this paper:

\begin{proof}{\bf of Theorem \ref{thm:convind}} We can write:
\begin{eqnarray*}
D(p_{n,k,t} \| \gibb{k} ) & = &
\int \frac{ q_{n,k}(\vc{y},nt)}{w_n(nt)}
\rat{q_{n,k}(\vc{y},nt)}{\gibb{k}(\vc{y}) w_{n-k}(nt-R_k(\vc{y}))} 
d\vc{y} \\
& &  +  \int p_{n,k,t}(\vc{y}) \rat{w_{n-k}(nt-R_k(\vc{y}))}{w_n(nt)} 
d\vc{y}.
\end{eqnarray*}
Proposition \ref{prop:logsum} below deals with the first term.
Lemmas \ref{lem:meanbd} and \ref{lem:intbd} imply Proposition 
\ref{prop:df}, providing a bound uniformly in $z$ on 
$\log( w_{n-k}(z)/w_n(nt))$, so we can deal with the second term.
\end{proof}
An analogue of the discrete log-sum inequality holds for integrals, and hence 
we deduce a projection inequality: 
\begin{proposition} \label{prop:logsum} 
Under the previous definitions:
$$ \int \frac{ q_{n,k}(\vc{y},nt) }{w_n(nt)}
\rat{ q_{n,k}(\vc{y},nt)}{\gibb{k}(\vc{y}) w_{n-k}(nt-R_k(\vc{y}))} d\vc{y}
\leq  D(p \| h_n).$$
\end{proposition}
\begin{proof} 
Given integrable functions $g(x),h(x)$, normalising to get probability 
densities
$p(x) = g(x) / \int g(x) dx$, $q(x) = h(x) / \int h(x) dx$, the Gibbs 
inequality gives:
$$ 0 \leq D( p \| q) \int g(x) dx = 
 \int g(x) \log \left( \frac{g(x)}{h(x)} \right) dx -  \left(\int g(x) dx 
\right) \log \left( \frac{\int g(x) dx}{\int h(x) dx} \right).  $$
Now, writing $\vc{x} = (\vc{y}, \vc{u})$,
where $\vc{y} \in \re^k$, $\vc{u} \in \re^{n-k}$, notice that:
\begin{eqnarray*}
 q_{n,k}(\vc{y},nt) & = & \int (\theta p)(\vc{y},\vc{u}) \I(R_k(\vc{y})
+ R_{n-k}(\vc{u}) = nt) d\vc{u}\\
 \gibb{k} (\vc{y}) w_{n-k}(nt-R_k(\vc{y})) 
& = & \int (\theta h_{n,t})(\vc{y},\vc{u}) \I(R_k(\vc{y})
+ R_{n-k}(\vc{u}) = nt) d\vc{u}
\end{eqnarray*}
Hence we deduce that for each $\vc{y}$:
\begin{eqnarray*}
\lefteqn{ \frac{ q_{n,k}(\vc{y},nt) }{w_n(nt)} 
\rat{ q_{n,k}(\vc{y},nt)}{\gibb{k}(\vc{y}) w_{n-k}(nt-R_k(\vc{y}))}} \\
& \leq &  \int \frac{ (\theta p)(\vc{y},\vc{u})}{w_n(nt) } 
\rat{(\theta p)(\vc{y},\vc{u})}{(\theta h_{n,t})(\vc{y},\vc{u})}
\I(R_k(\vc{y}) + R_{n-k}(\vc{u}) = nt) d\vc{u}\\
& = & \int p(S_n(\vc{y},\vc{u})) 
\rat{p(S_n(\vc{y},\vc{u}))}{h_{n,t}(S_n(\vc{y},\vc{u}))} 
\I(R_k(\vc{y}) + R_{n-k}(\vc{u}) = nt) d\vc{u}
\end{eqnarray*}
Integrating with respect to $\vc{y}$, we obtain:
\begin{eqnarray*}
\lefteqn{ \int \frac{ q_{n,k}(\vc{y},nt) }{w_n(nt)} 
\rat{ q_{n,k}(\vc{y},nt)}{\gibb{k}(\vc{y}) w_{n-k}(nt-R_k(\vc{y}))} d\vc{y}} 
\\
& = & \int p(S_n(\vc{y},\vc{u})) \rat{p(S_n(\vc{y},\vc{u}))}
{h_{n,t}(S_n(\vc{y},\vc{u}))} 
\I(R_k(\vc{y}) + R_{n-k}(\vc{u}) = nt) d\vc{u} d\vc{y} \\
& = & \int p(S_n(\vc{x})) \rat{p(S_n(\vc{x}))}{h_{n,t}(S_n(\vc{x}))} 
\I(R_n(\vc{x}) = nt) d\vc{x} \\
& = &  D(p \| h_n),
\end{eqnarray*}
as required.\end{proof}
\begin{corollary}\label{cor:unif}
Assume $f \in \mathcal{F}$ and taking $t = \ene(\gibb{1})$, given the 
density $p = h_{n,t}$  on 
${\mathcal{S}}_n(t)$, then the projection $p_{n,k,t}$ satisfies:
$$ D(p_{n,k,t} \| \gibb{k} ) \leq 
\rat{n}{n-k} + \frac{2}{\sqrt{n}/C  - 1}.$$ \end{corollary}
Observe that $d_{TV}(f,g) = \int | f(x) - g(x) | dx \leq 
\sqrt{2 D(f \| g)}$ (see \cite{kullback}), and hence in this case
then the rate of convergence in total variation distance is $O(1/\sqrt{n})$,
as opposed to the $O(1/n)$ which Diaconis and Freedman establish.
This difference can be attributed to the fact that we approximate the
densities $w_n$, rather than being able to obtain exact bounds on them.
\section{A converse}
Given the density $p= h_{n,t}$ on ${\mathcal{S}}_{n}(t)$
we can see from Corollary \ref{cor:unif} that as 
$n \tends$, and $k/n \rightarrow 0$ then $D(p_{n,k,t} \| \gibb{k}) 
\rightarrow 0$.  This is also Diaconis and Freedman's \cite{diaconis} 
necessary and sufficient condition for convergence in total variation distance 
in the spherical case. We show that
this condition holds for more surfaces than that:
\begin{proposition} Given $f \in \mathcal{F}$,
consider the uniform probability density $h_{n,t}$ on ${\mathcal{S}}_{n}(t)$,
and consider $p_{n,k,t}$, the distribution of its first $k$ coordinates. If
$p_{n,k,t} \rightarrow \gibb{k}$ in total variation distance
then $k/n \rightarrow 0$.
Hence the stronger result of $D(p_{n,k,t} \| \gibb{k}) \rightarrow 0$ 
implies that $k/n \rightarrow 0$.
\end{proposition}
\begin{proof}
Now considering $\vc{Y}$ with density $\gib$, by conditioning and 
independence, the density $r_k(y)$ of  $(R_k(\pi_{n,k} \vc{Y}) | 
R_n(\vc{Y})=n t)$ satisfies:
$$ r_k(y) = \frac{w_k(y) w_{n-k}(nt -y)}{w_n(n t)}.$$
Hence since total
variation distance is reduced by projection,  we deduce that for any set $L$:
$$ d_{TV}(p_{n,k,t}, \gibb{k}) \geq d_{TV}( r_k , w_k)
\geq 2 \int_L w_k(y) \left( \frac{w_{n-k}(nt -y)}{w_n(n t)} -1 \right) 
dy.$$ 
Now, using the estimates of the previous section, we know that choosing the
interval $L = (k t - \epsilon \sqrt{n-k}, k t + \epsilon \sqrt{n-k})$
will ensure that the term in brackets is close to $\sqrt{n/n-k}-1$.
Furthermore for $n \gg k$, $w_k(L)$ will be close to 1, so we deduce the 
result.
 \end{proof} 
\section*{Acknowledgements}
I wish to thank Yuri Suhov and Geoffrey Grimmett for advice, 
assistance and support. This work was supported by an EPSRC grant, 
Award Reference Number 96000883, and EC Grant `Training 
Mobility and Research' (Contract ERBMRXT-CT 960075A). I would also like to
thank Professor Hans-Otto Georgii of LMU, M\"{u}nchen for useful discussions.
Professors Joel Lebowitz and Eric Carlen also provided extremely helpful 
advice.

\end{document}